\documentclass[12pt]{article}
\usepackage{amssymb}
\input epsf

\def \be{\begin{equation}}
\def \ee{\end{equation}}
\def \berr{\begin{eqnarray}}
\def \err{\end{eqnarray}}
\def \bea{\begin{array}}
\def \ea{\end{array}}
\def \nn{\nonumber}
\def \a{\alpha}

\def \d{\delta}

\def \l{\lambda}

\def \dl{\partial}

\def \vareps{\varepsilon}
\def \A{{\cal A}}

\def \H{{\cal H}}

\def \U{{\cal U}}

\def \({\left(}
\def \){\right)}
\def \<{\langle}
\def \>{\rangle}
\def \[{\left[}
\def \]{\right]}

\def \obar{\overline}

\def \del{\partial}
\def \tens{\mathop{\otimes}}
\def \btens{\tens_h}

\def\mf#1{{\mathbb #1}}

\def\R{{\mf{R}}}
\def\C{{\mf{C}}}

\def\mg{\mathfrak{g}}
\newcommand \one{{\bf 1}}

\def\c#1{{\cal #1}}

\def\smash{\mbox{$\,\rule{0.3pt}{1.1ex}\!\times\,$}}

\newcommand{\sect}[1]{\setcounter{equation}{0}\section{#1}}

\evensidemargin \oddsidemargin
\begin{document}

\begin{titlepage}

\begin{center}  

\Large\bf
Propagator on the $h$--deformed Lobachevsky plane
 \\[4ex]
\normalsize \rm
J.\  Madore$^a$\footnote{John.Madore@th.u-psud.fr}
and H. \ Steinacker$^b$\footnote{Harold.Steinacker@physik.uni-muenchen.de} 
 \\[2ex] 
{\small\it 
        ${}^a$Max-Planck-Institut f\"ur Physik\\
        (Werner-Heisenberg-Institut)\\
        F\"ohringer Ring 6, D-80805 M\"unchen  \\
        and \\
        Laboratoire de Physique Th\'eorique et Hautes Energies\\
        Universit\'e de Paris-Sud, B\^atiment 211, F-91405 Orsay \\[2ex]

        ${}^b$Sektion Physik der Ludwig--Maximilians--Universit\"at M\"unchen\\
        Theresienstr.\ 37, D-80333 M\"unchen  \\[1ex] }
\end{center}

\vspace{5ex}

\begin{abstract}
The action of the isometry algebra $U_h(sl(2))$ on the $h$--deformed 
Lobachevsky plane is found. The invariant distance and 
the invariant $2$--point functions are shown
to agree precisely with the classical ones. The propagator 
of the Laplacian is calculated  explicitely. It is invariant only 
after adding a ``non--classical'' sector to the Hilbert space.
\end{abstract}

\parskip 4pt plus2pt minus2pt
\vfill
\noindent
PACS classification 11.10.Kk, 03.65.Fd 

\end{titlepage}

\sect{Introduction}

The $h$--deformed Lobachevsky plane was introduced by Demidov et al.
\cite{DemManMukZhd90} and by Manin \cite{Man91}. Its function algebra
is covariant under $Sl_h(2,\mf{R})$, which is a triangular 
Hopf algebra, sometimes called the Jordanian
deformation of $Sl(2,\mf{R})$. As opposed to the
$q$--deformed quantum groups, it is triangular, which means
that the deformation from the classical case 
is less severe. In fact it is  known that $U_h(sl(2))$ 
is related to its undeformed counterpart by a twist \cite{abdel,kulish}.
While this might suggest that the deformation is almost trivial 
in some sense, the problem of defining suitable spaces of functions,
in particular Hilbert spaces, which can be relevant to physical systems
is not at all trivial.
In fact it will turn out that in order to find an
invariant propagator, a certain ``non--classical'' sector must be added 
to the Hilbert space, which disappears in the classical limit.

The main goal of this paper is to calculate explicitly the 
invariant propagator on the $h$--deformed Lobachewski plane. 
The first observation is that the well--known 
covariance algebra $U_h(sl(2))$ does  not preserve
the metric structure, but only the symplectic structure. Therefore
in Section 2, we first
determine the 3--parameter ``group'' of isometries, which turns out
to be again $U_h(sl(2))$, but with a different action on the space
which corresponds to the well-known
fractional transformations of the upper half--plane.

In order to define $n$--point ``functions'' in a covariant way, 
in Section 3 we introduce 
braided copies of the $h$--deformed Lobachevsky plane. This allows
to determine invariant functions, and in particular the
invariant distance between 2 ``points''. The distance 
turns out to involve only a commutative subalgebra of the complete
algebra, and  agrees precisely with the classical one. 

In Section 5, we calculate the propagator of the $h$--deformed 
Laplacian explicitly. When based on a naive generalization of the 
Hilbert space of modes of the undeformed case, it turns out 
that the propagator is not invariant under $U_h(sl(2))$. 
It does become invariant only after adding another, 
``nonclassical''  sector to the Hilbert space. 
This situation is reminiscent
of a similar phenomenon on  $q$--deformed quantum spaces \cite{schwenk}
and shows that the $h$--deformation is not quite a trivial one.
The propagator on the extended Hilbert space 
then turns out to agree formally with the classical one.

This result should be compared with that of a recent work 
\cite{madore_etc}  where
the propagator on the $h$--deformed Lobachevsky plane has been found
to be finite provided one uses the {\em undeformed} tensor
product. Of course, this breaks the covariance under $U_h(sl(2))$. 
In our covariant treatment, the propagator turns out
not to be regularized. This means that either $h$--deformation is 
not strong enough to regularize the UV divergencies, or that the 
different copies of the Hilbert space should {\em not} be implemented
via  the braided tensor product.

\sect{The isometries of the $h$--deformed Lobachevsky plane}

The $h$--deformed  Lobachevsky plane \cite{Agh93} can be  defined 
\cite{cho_madore} to be the formal
$*$-algebra $\A$ generated by two hermitian elements $x$ and $y$
which satisfy the commutation relation 
\be
[x,y] = - 2 i h y                                       \label{lob_alg}
\ee
where $h \in \mf{R}$ and the factor $-2$ is present for historical
reasons. We shall suppose that $h > 0$. 
Throughout this paper, a ``function'' on the $h$--deformed  
Lobachevsky plane is understood to be an element of 
$\A$ or a suitable completion thereof.

Using the variable $z=x+iy$, this becomes
\be
\[z,\bar{z}\] = 2ih(z-\bar{z}) .                         \label{z_CR}
\ee
Introducing variables $r,s$ by $x=rs^{-1}+\frac 12 ih$ and 
$y=s^{-2}$, the above commutation relation becomes
\be
\[r,s\] = ih s^2.
\label{rs_alg}
\ee
In terms of these variables, it is easy to check that the algebra
is covariant under $Sl_h(2,\mf{R})$, i.e. there is a coaction 
$\Delta:\A \rightarrow {\rm Fun}(Sl_h(2,\mf{R})) \tens \A$ given by
\be
\(\bea{c}r\\s \ea\) \rightarrow \(\bea{cc} A & B \\ C & D \ea \) 
  \stackrel{\cdot}{\tens}  \(\bea{c}r\\s \ea\)      \label{rs_coaction}
\ee
where the algebra ${\rm Fun}(Sl_h(2,\mf{R}))$ is the 
$h$--deformed (Hopf) $*$--algebra
of functions on $Sl(2,\mf{R})$ generated by the hermitian elements $A,B,C,D,$ 
with relations
\berr
\[A,B\] & = & ih\delta - ih A^2, \nonumber \\
\[A,C\] & = &  ih C^2,          \nonumber \\
\[A,D\] & = & ih CD -ih CA,      \nonumber \\
\[B,C\] & = & ih CD + ih AC,     \nonumber \\
\[B,D\] & = & ih D^2 - ih \delta,\nonumber \\
\[C,D\] & = & - ih C^2,         \nonumber
\err
where the quantum determinant
\be
\delta   = AD - CB - ih CD = DA - CB - ih CA
\ee
is central and set equal to one.

The $R$-matrix associated to this quantum group, which solves
the quantum Yang-Baxter equation 
\be
\hat{R}_{12}\hat{R}_{23}\hat{R}_{12}
= \hat{R}_{23}\hat{R}_{12}\hat{R}_{23},                        
\ee
is given by
\be
\hat{R} = \left( \begin{array}{cccc}
             1  & -ih & ih & -h^2 \\
             0  &  0  & 1 & ih \\
             0  &  1 & 0 & -ih  \\
             0  &  0 & 0 & 1
                \end{array}
              \right).                                   
\ee
It is triangular, i.e. $\hat{R}^2 = 1$, which also holds for 
the higher representations.
The associated calculus and Laplacian have been worked out 
elsewhere \cite{Agh93,cho,cho_madore}. We will thus be brief here.

The covariant differential calculus $(\Omega^*(\c{A}_h),d)$
on a quantum plane can be
found~\cite{Agh93} by the method of Wess and Zumino~\cite{WesZum90}.
For $r^i = (r, s)$ and $\xi^i = dr^i = (\xi, \eta )$ we have
\berr
&r^a r^b =        \hat{R}^{ab}{}_{cd}   r^c  r^d, \hspace{1cm}
&r^a \xi^b =      \hat{R}^{ab}{}_{cd} \xi^c  r^d,  \nonumber   \\
&\xi^a \xi^b =  - \hat{R}^{ab}{}_{cd} \xi^c\xi^d, \hspace{1cm} \label{calc}
&\partial_a x^b = \delta_a^{b}
                + \hat{R}^{bd}{}_{ac}x^c \partial_{d}.        
\err
Explicitely, the second and third equations are 
\be
\begin{array}{ll}
\[r, \xi\] = - ih \xi s +ih \eta r - h^2 \eta s, & \[r, \eta\] = ih \eta s, \\
\[s, \xi\] = - ih \eta s, &\[ s, \eta\] = 0,
\end{array}                                                  
\ee
and
\be
\xi^2 =  ih \xi \eta, \qquad \xi \eta = - \eta \xi,  \qquad
\eta^2 = 0.                                                  
\ee
One can also introduce 
\cite{cho_madore} a frame or Stehbein $\theta^a$ defined by
\be
\theta^1 = y^{-1} dx, \qquad \theta^2 = y^{-1} dy. \qquad    \label{theta}
\ee 
They satisfy the commutation relations
\be
f \theta^a = \theta^a f, \qquad f \in \c{A}_h                  \label{2.2}
\ee
as well as the quadratic relations
\be
(\theta^1)^2 = 0, \qquad (\theta^2)^2 = 0, \qquad
\theta^1 \theta^2 + \theta^2 \theta^1 = 0.
\ee
More details of this have been given elsewhere~\cite{cho,cho_madore}.
One can also define a Hodge--star operator as
\be
\ast(\theta^1) = \theta^2, \; \ast(1) = \theta^1 \theta^2;
\label{hodge}
\ee
we will see in section 3 that this is also the correct covariant 
definition.

By construction, the coaction of ${\rm Fun}_h(Sl(2,\mf{R}))$ on $\A$ 
preserves the symplectic structure; however it is not the 
group of isometries. This is so even classically, as was already noted in 
\cite{cho_madore}. To find the correct isometries, consider first the 
commutative limit, where the metric is $ds^2 = y^{-2} (dx^2 + dy^2)$.
The isometries are the well--known fractional transformations
$$
z \rightarrow z' = \frac{Az+B}{Cz+D}\qquad {\rm with} \quad 
\(\bea{cc}A & B\\C & D\ea\) 
\in {\rm Fun}(Sl(2,\mf{R})),
$$
where $z=x+iy$.
In this section, we will find a similar
transformation in the noncommutative case, such that the
isometries have the structure of $SL_h(2,\mf{R})$, but with a 
different (co)action than the symplectic one (\ref{rs_coaction}).  
The metric will turn out to be the same as in the commutative case.

Since we later wish to determine the functions which are invariant under
the isometries, it is more useful to have an action of the
universal enveloping algebra $U_h(sl(2))$ on $\A$ rather than a
coaction of ${\rm Fun}_h(Sl(2))$. Since these are dually paired Hopf algebras,
a left (respectively right)  coaction of ${\rm Fun}_h(Sl(2))$
corresponds to a right (respectively left) action of $U_h(sl(2))$, 
see for example \cite{majid_JMP}. The resulting cross--product
algebra is given in (\ref{cross_alg}); we take a small detour and 
explain the steps leading to this algebra.

To find this dual action, we look for
variables such the the above fractional transformation of $z$
becomes linear.
First, we introduce different variables
$z_1, z_2$ for $\A$ which satisfy 
\be
\[z_1,z_2\] = 2ih (z_1- z_2),
\label{z_1z_2}
\ee
with star--structure  $z_1^* = \obar{z_1} =  z_2 + ih$;
we use a bar to denote the star of $z$.
Then (\ref{z_CR}) is 
recovered for $z = z_1 + \frac{ih}2$.
One can easily check (and it will become evident below)
that this is consistent with the 
following coaction of ${\rm Fun}_h(Sl(2))$:
\be
z_i \rightarrow (Az_i+B)(Cz_i+D)^{-1}                    \label{z_coaction}
\ee
for $i=1,2$. A similar coaction 
for certain $q$--deformed algebras has been considered in \cite{pmho}. 
While this form is very appealing, it is somewhat formal, and
it is not immediately clear how to
translate it into an action of $U_h(sl(2))$ which will needed below.
To find this, we introduce yet another set of (auxiliary) generators. 
Consider $u_i,v_i$ with $\[u_i,v_i\] = ih v_i^2$
for $i=1,2$, which are
covariant under the linear coaction of ${\rm Fun}_h(Sl(2))$ as in
(\ref{rs_alg}) and (\ref{rs_coaction}),
and let $z_i=u_i v_i^{-1}$.
Furthermore, we 
impose the commutation relations
\berr
[u_1, u_2] &=& ih u_1 v_2 -ih v_1 u_2 +h^2 v_1 v_2, 
\nonumber \\ \ 
[u_1, v_2] &=& ih v_1 v_2, \nonumber \\ \ 
[v_1, u_2] &=& -ih v_1 v_2, \nonumber \\ \ 
[v_2,v_2]  &=& 0.                               \label{basic_braid}
\err
They are consistent with the linear coaction of ${\rm Fun}_h(Sl(2))$ as
will be explained in the next section, and imply
(\ref{z_1z_2}). The star structure $u_1^* = u_2$, 
$v_1^* = v_2$ implies $z_1^* = \obar{z_1} = z_2 +ih$ as above, and the linear
coaction on $u_i,v_i$ obviously induces (\ref{z_coaction}).
In this linear form, we can find the dual action of $U_h(sl(2))$,
and then restrict it to the original algebra generated by $z, \obar{z}$.

We recall the definition of $U_h(sl(2))$:
It is the Hopf algebra with generators $\{J^{\pm}, J^3\}$ and relations 
 \cite{ohn} 
\berr
\[J^3,J^+\] &=& 2h^{-1} \sin(hJ^+), \nonumber\\ 
\vspace{2pt}
\[J^3,J^-\] &=& -\big[\cos (hJ^+)J^- + J^-\cos (hJ^+)\big],  \nonumber\\
\vspace{2pt}
\[J^+,J^-\] &=& J^3                                     \label{U_h_alg}
\err
and
\berr
\Delta J^+ &=& J^+\otimes 1+1\otimes J^+,  \qquad
           \Delta J^j=J^j\otimes e^{-ihJ^+}+ e^{ihJ^+}\otimes J^j, \nonumber\\
\vspace{4pt}
\vareps(X) &=& 0, \quad  S(X)=-e^{-ihJ^+}Xe^{ihJ^+},  
\err
where $j \in \{ -,3\}$ and $X\in\{ J^+,J^-,J^3\}$.
This is obtained from the result of Ohn \cite{ohn} by the 
replacement $h \rightarrow -ih$.
It is a $\ast$--Hopf algebra  with the reality structure
\be
(J^{\pm})^{\ast} = -J^{\pm}, \quad (J^3)^{\ast} = - J^3,
\label{h_star}
\ee
which defines $U_h(sl(2,\mf{R}))$.  Introducing 
\be
G=e^{-ihJ^+},
\ee
this becomes
\berr
[G,J^3]   &=& 1-G^2, \nonumber \\  \ 
[J^3,J^-] &=& -\frac 12 \big[(G+G^{-1})J^- + J^-(G+G^{-1}) \big]\nonumber\\ \  
[G,J^-]   &=& -\frac{ih}2(G J^3 + J^3 G), 
\err
and
\berr
\Delta G &=& G \tens G,  \qquad \Delta J^j=J^j\otimes G+
                         G^{-1} \otimes J^j,                  \nonumber\\
\vspace{4pt}
\vareps(X) &=& 0,  \qquad S(X)=-G X G^{-1},                  
\err
where $j \in \{-,3\}$ and $X \in\{ J^+,J^-,J^3\}$.  
Given a (left or right) action of $U_h(sl(2))$ on $\A$, one can 
always define a
cross--product algebra $U_h(sl(2)) \smash \A$. As a vector space, 
this  is  $U_h(sl(2)) \tens \A$,
equipped with an algebra structure defined by
$ua = (u_{(1)}\cdot a)u_{(2)}$, where the dot denotes the left action of
$u \in U_h(sl(2))$ on a representation; similarly for a right action.
Here $u_{(1)} \tens u_{(2)}$ denotes the coproduct of $u$.
Conversely, the  left action $u \cdot a$ can be extracted by commuting
$u$ to the right and then applying the counit $\vareps$ of $U_h(sl(2))$ 
from the right.

The dual pairing of $U_h(sl(2))$ with ${\rm Fun}_h(Sl(2))$ has been
given in \cite{dobrev}. Using this, 
it is easy to find the dual action of $U_h(sl(2))$ 
on the variables $u_1,v_1, u_2, v_2$. 
This defines a cross--product algebra as explained above, 
which can be expressed in terms of the original variable $z$.
The resulting algebra is
\berr
[z,J^+]   &=&  -1,                \nonumber \\ \ 
[z,J^- G] &=& (z^2 - ihz + \frac{h^2}4) G^2,      \nonumber\\  \ 
[z,J^3 G] &=&  (2z-ih) G^2                        \label{cross_alg}
\err
and the same relations for $\obar{z}$. 
One can check explicitely that this algebra is consistent with the relations
(\ref{z_CR}), the relations of $U_h(sl(2))$, 
and also with the star structure $z^{\ast} = \obar{z}$. 
Therefore it is a consistent cross--product algebra for the 
$h$--deformed Lobachevsky plane.
In terms of $x$ and $y$ with $z=x+iy, \obar{z} = x-iy$, it becomes
\berr
[x,J^+]   &=& -1, \qquad[y,J^+]  = 0,  \nonumber \\  \ 
[x,J^-G]  &=& (x^2 - y^2 - ihx + \frac{h^2}4) G^2, \nonumber \\  \ 
[y,J^-G]  &=& (2xy + ihy) G^2 \nonumber \\ \ 
[x,J^3G]  &=& (2x-ih) G^2,  \nonumber \\ \ 
[y,J^3G]  &=& 2y G^2 .                          \label{cross_alg_real}
\err
With 
\be
D_x f(x) :=\frac{f(x) - f(x-2ih)}{2ih},
\ee
one finds for functions (power series) in $x$ and $y$ 
\berr
 [f(x), J^+]  &=& -\del_x f(x), \qquad \qquad \qquad  [f(y),J^+] = 0, \nonumber\\ \ 
[f(x),J^3G] &=& (2x-ih) D_x f(x) G^2 , \quad
   [f(y),J^3G] = 2y \del_y f(y) G^2   \nonumber\\ {}
[f(x),J^-G] &=& ((x^2-ihx + \frac{h^2}4) D_x f(x) - D_x f(x+2ih)y^2) G^2, 
                                                            \nonumber \\ {}
[f(y),J^-G] &=& \(2xy\del_y f(y)+ih\(2(y\del_y)^2 - y\del_y\)f(y)\) G^2  
\label{CR_functions}
\err
and 
$$
Gf(x) = f(x-ih)G.
$$  
We claim that this cross--product algebra implements the 
isometry algebra $U_h(sl_h(2))$ on the $h$--deformed Lobachevsky plane.
In the limit $h \rightarrow 0$, the generators $J^+,-J^3G,$ and $ -J^-G$ 
clearly become 
the classical generators $\del_x, 2x \del_x + 2y \del_y,$ and 
$(x^2 - y^2) \del_x + 2xy \del_y$ of the algebra $sl(2,\R)$ of isometries.
Moreover in Section 3.1, we shall find an explicit 
expression for the $h$--deformed distance,
which is invariant under the action (\ref{cross_alg_real}) of $U_h(sl_h(2))$.
The necessary tools will be provided in the next section. 

\sect{Braided copies of the $h$--deformed Lobachevsky plane}

In order to write down functions of several variables such as 
$n$--point functions, 
one should introduce several copies of the algebra of functions, and
combine them into a bigger algebra. 
Recall the  classical case: if $\A$ is a representation of some
Lie algebra $\mg$ compatible with the algebra structure of $\A$, 
more precisely $\A$ is a $\mg$ --module algebra, this is
easy to do: define $\A^{\tens} :=\A \tens ... \tens \A$, and let
$\A^{(n)}:=1 \tens ... \tens \A \tens 1 \tens ...$ be the $n$th copy of 
$\A$. Then $\A^{\tens}$
is naturally an algebra (the tensor product algebra)  
by component--wise multiplication, and
$\A^{(n)}$ commutes  with $\A^{(m)}$ if $n \neq m$. $\A^{\tens}$
carries the tensor product representation of $\mg$, and 
its algebra structure is compatible with this representation. 

If  $\A$ is covariant under a Hopf algebra $\U$
which is not cocommutative, this standard algebra stucture on $\A^{\tens}$
is not compatible with the coaction. 
However if the Hopf algebra is quasitriangular
with universal $R$--``matrix''
${\cal R} = {\cal R}_1 \tens {\cal R}_2 \in \U \tens \U$ (in short--hand
notation), then there is a 
standard way to define a modified algebra structure on $\A^{\tens}$,
the so--called ``braided tensor product'' \cite{majid}:
it is defined by
$a^{(n)} a^{(m)} := 1 \tens ... \tens a^{(n)} \tens 1 \tens ...\tens a^{(m)}
\tens 1 ...$ if $n \leq m$, and
$a^{(n)} a^{(m)} := 1 \tens ... \tens {\cal R}_1 \cdot a^{(m)} \tens 1 
\tens ...\tens {\cal R}_2 \cdot a^{(n)}
\tens 1 ...$ if $n > m$, where $a^{(n)} \in \A^{(n)}$.
This is compatible with the 
action of the quasitriangular Hopf algebra $\U$;
to avoid confusion, we will denote it by $\A^{\btens}$.

Since $U_h(sl(2))$ is in fact triangular, i.e.
${\cal R}_{12} {\cal R}_{21} = \one$, the above
definition can be written as a commutation relation
$a^{(n)} a^{(m)} = ({\cal R}_1 \cdot a^{(m)}) ({\cal R}_2 \cdot a^{(n)})$
in $\A^{\btens}$ whenever $n \neq m$. This is a considerable
simplification over the quasitriangular case, where one has to 
distinguish between $n>m$ and $n<m$.
Notice that the commutation relations between functions and the generators
of forms in the calculus (\ref{calc}) are precisely of this kind.

For the $h$--deformed Lobachevsky plane, we have 2 different actions
of $U_h(sl(2))$ acting on it, the symplectic one dual to (\ref{rs_coaction})
and the (tentative) isometries corresponding to (\ref{cross_alg}) 
found in the previous section. Thus it is not clear a priori how to 
proceeed. What we will do is to define first the braided 
algebra $\A^{\btens}$ 
for the symplectic action since that one is much simpler, and
verify that it in fact compatible with the action of the 
isometries as well; this is not obvious a priori.

For simplicity, introduce just 2 copies of $\A$, i.e. $x = x\tens 1$ and
$x' = 1 \tens x$.
In terms of the variables $r,s$ of section 2 
with $x=rs^{-1}+\frac 12 ih$ and 
$y=s^{-2}$, this definition leads to
\berr
[r, r'] &=& ih r s' -ih s r' + h^2 s' s,  \nonumber \\ \ 
[r, s'] &=& ih s' s, \nonumber \\ \ 
[s, r'] &=& -ih s' s, \nonumber \\ \ 
[s,s']  &=& 0, \nonumber                             
\err
which is the same as (\ref{basic_braid}).
In terms of $x$ and $y$, this becomes
\berr
\[x,x'\]  &=& 2ih x -2ih x' \nonumber\\ \ 
\[x,y'\]  &=&  -2ih y'      \nonumber\\ \ 
\[y,x'\]  &=& 2ihy          \nonumber\\ \ 
\[y,y'\]  &=& 0.                         \label{braid_alg}
\err
It is somewhat disturbing that the commutator of $x$ and $x'$ does not
vanish as $x-x'$ becomes large, but this is required by covariance; 
we will come back to that later. In the complex variables
$z=x+iy, \quad z' = x' + iy'$, one obtains
\be
\[z,z'\] = 2ih (z-z'),
\ee
which explains the relation (\ref{z_1z_2}).
This algebra is by construction consistent with the coaction
of the ``symplectic'' $Fun(Sl_h(2))$ (\ref{rs_coaction}) respectively 
its dual. 
It can now be checked by a lenghty but straightforward calculation 
that these relations are also compatible 
with (\ref{cross_alg}), extended to both copies $z$ and $z'$.
This is not obvious a priori.
A somewhat similar observation has been made \cite{pmho} 
for the $q$--deformed case 
in terms of the fractional transformations considered in section 2.

The concept of braided copies of a covariant algebra is also relevant
if one tries to define a Fock space of creation-- and anihilation 
operators which are covariant under some quantum group. 
In general, it is not obvious then how to define a totally
symmetric or antisymmetric Hilbert space, since the deformed analogue of the
permutation operator, $\hat{R}$, has eigenvalues which are different
from $\pm 1$. In the triangular case, this problem does not
occur, since $\hat{R}^2 = \one$ by definition. From this point of view,
the triangular case seems particularly well suited  to formulate 
a Quantum Field Theory. Perhaps however  
triangular Hopf algebras are not  a ``sufficiently'' nontrivial
deformation in order to improve the UV behaviour of the commutative limit. 
To obtain some insight into this question was one of the motivations of
the present work. 

\subsection{Invariant distance}

To make the algebra (\ref{braid_alg}) more transparent, we define 
$$
\d x = x-x', \quad  \d y = y-y'
$$
and
$$
\obar{x} = \frac 12 (x+x'), \quad \obar{y} = \frac 12 (y+y').
$$
In terms of theses variables, one finds as only nonzero commutators
\berr
[\obar{x}, \obar{y}] &=& -2ih \obar{y}, \nonumber \\ \ 
[\obar{x}, \d x]     &=& -2ih \d x, \nonumber \\ \ 
[\obar{x}, \d y]     &=& -2ih \d y.  
\err
Notice that these are the same relations as for the calculus where $\d x, \d y$
are replaced by $dx$ and $dy$. In particular, 
$\obar{y}^{-1} \d x$ and $\obar{y}^{-1} \d y$ play the role of
the Stehbein (\ref{theta}), and they commute with $\obar{x}$ and $\obar{y}$.
Thus there is only one nontrivial commutator among the four generators
of $\A \btens \A$ as opposed to the case $\A \tens \A$, where the
propagator is regularized \cite{madore_etc}.
In particular, it is somewhat counterintuitive that 
the ``relative'' and ``average'' coordinates
do {\em not} mutually commute (cp. \cite{madore_etc}); 
again, this is forced upon us by the covariance requirement.

The geodesic distance of 2 points $(x,y)$ and $(x',y')$ on the 
classical Lobachavsky plane is given by \cite{Kub88}
\be
d = \cosh^{-1} \(1+ \frac 1{2yy'}((\d x)^2 + (\d y)^2) \).       
        \label{distance}
\ee
The subalgebra generated by $\obar{y}, \d x, \d y$
is abelian in the $h$--deformed case as well, and
it can be checked that the same expression  is also
invariant under the $h$--deformed isometries. In fact,
$$
y^{-1} y'^{-1} ((\d x)^2 + (\d y)^2) 
$$
commutes with $U_h(sl(2))$ in the cross--product algebra (\ref{cross_alg}),
and together with 1 generates the center of $\A{\btens} \A$.
Again, this is more easily seen in terms of the 
fractional transformation of Section 2, but less rigorous.
Thus we shall define (\ref{distance}), which is
an invariant 2--point function, 
to be the invariant distance function on the
$h$--deformed Lobachevsky plane. 

Of course the same considerations apply for the case of several variables.
It is then clear that the set of invariant $n$--point functions is
the same as  classically, i.e. they are precisely the functions 
which depend only on the relative distances of any pairs of
variables, defined as in (\ref{distance}). 

One can check that the Hodge--star (\ref{hodge}) is invariant
under the isometries as well.

\sect{Invariant functionals and inner products}

In quantum mechanics, symmetries are implemented as unitary
transformations on a Hilbert space. 
To realize this and similarly to define invariant propagators,
one has to find a positive definite inner product on a vector space
which is invariant under the symmetry.
Such invariant inner products are naturally obtained 
from a positive state, which should satisfy
\be
\< u \cdot a \> = \vareps(u) \< a \>        \label{invar_state}
\ee
and $\<a\>^* = \<a^*\>$, where $a$ is an element of a star algebra
$\A$ and $u$ an element of a symmetry (Hopf) algebra $\U$.
Since we are considering spaces of functions on a (noncommutative) 
manifold, this can be considered as an invariant functional.
Of course, $\<\;\>$ is defined only on a certain subset of 
``measurable'' elements of a suitable completion of 
$\A$, as classically. This will become clear in our example.

It is useful to formulate this within the framework of cross--product 
algebras. Any state (functional) on $\A$  
defines a state on $U_h(sl(2)) \smash \A$, written
as $\< ua \>$, by $\< ua \> := \epsilon(u) \<a\>$.
If the state on $\A$ is invariant, this implies 
using some standard identities of Hopf algebras that 
$\<uav\> = \epsilon(u)\epsilon(v) \<a\>$
for any $u,v \in U_h(sl(2))$, and in particular
\be
\< [u, a] \> =0                                 \label{invar_CR}
\ee
for any $u \in U_h(sl(2))$.
Conversely, the state $\<\;\>$ on $\A$ is invariant if (\ref{invar_CR}) holds.
The latter form is quite intuitive, and well suited for our situation.
We will work with this formalism from now on. 

As usual, each invariant state induces an invariant inner product 
as follows:
\be
\<f,g\> = \<f^* g\>.                         \label{inner_prod}
\ee
It is invariant, because 
$\<f,u \cdot g\> = \<f^*u g\> = \<(u^*\cdot f)^* g\>$,
using $u\cdot f = u_1 f S u_2$ and standard identities of 
Hopf--$\ast$  algebras.

The conditions (\ref{invar_CR}) for the subalgebra  of  $U_h(sl(2))$
generated by $G,J^-,J^3$ are 
\berr
\< [f,G] \>     &=& 0, \label{cond_1} \\
\< [f,J^3G]  \>     &=& 0, \label{cond_2} \\
\< [f, J^-G] \> &=& 0, \label{cond_3}
\err
for $f \in \A$. We will write the elements $f$  
in the form $f(x|y) = \sum f_{n,m} x^n y^m$, i.e. with $x$ to the left 
of $y$. After some calculations using (\ref{CR_functions}), they reduce 
to the following 2 requirements:
\be
\<f(x|y)\> = \<f(x+ih|y) \>                       \label{finite_transl}
\ee
and
\be
\< y \frac{\dl}{\dl y} f(x|y) \> = \<f(x|y) \>.
\label{cond_4}
\ee
Here $\frac{\dl}{\dl y} f(x|y)$ is the ordinary
differentiation w.r.t. $y$, after ordering the variables as above.
It turns out that (\ref{cond_3}) is a consequence of (\ref{cond_1}) and
(\ref{cond_2}). This is similar to the classical case, where the invariant
integral is also uniquely determined by 2 isometries, and automatically
respects the third.
The only difference to the classical case is that the translation
invariance w.r.t. $x$ is imposed only for a finite displacement rather
that for all. This of course comes from restricting ourselves to the 
algebra generated by $G, J^3,J^-$ rather than $J^+, J^3, J^-$, which
is consistent in the $h$--deformed case only.

One invariant functional satisfying theses conditions is now obvious:
It is simply the classical one. That is, consider the 
space $L^1(\mf{R}^2_+, d\mu)$
of functions $f(x,y)$ on the upper half plane 
which are integrable with respect to the measure
$d\mu = y^{-2} dx dy =\theta^1 \theta^2$. 
We write the functions (or more precisely a dense set of analytic
functions in $L^1(\mf{R}^2_+, d\mu)$) in the form $f(x|y)$ so 
that they define elements in (a completion of) $\A$, and we set
\be
\< f(x|y) \>_{(0)} :=\int f(x,y) d\mu.
\ee 
Invariance under $G$ follows by analytic continuation in $x$, e.g.
using the basis of Hermite--functions in $x$. 
In this way, we obtain a space
of functions on the $h$--deformed Lobachevsky plane 
which is isomorphic to $L^1(\mf{R}^2, d\mu)$. 
The corresponding Hilbert space 
will be defined  explicitely in the next section. 
At this point, the $h$--deformed case indeed
appears to be isomorphic to the undeformed case.

It will turn out, however, that this  ``classical'' Hilbert space is not
sufficient to obtain invariant propagators, but we shall be able to 
introduce ``extended'' Hilbert spaces by  
taking advantage of the weaker requirement (\ref{finite_transl}). 

Finally, one can define an integral of 2--forms $\a = f(x|y) \theta^1\theta^2$
by 
$$
\int \a = \<f(x|y)\>.
$$ 
It is easy to see that invariance of $\<\; \>$ 
under $G,J^-$ and $J^3$ is 
equivalent to Stokes theorem, $\int d\omega =0$, and that the adjoint of $d$ 
in the usual sense is indeed $\delta = \ast d \ast$. 


\sect{The propagator}

The $h$--deformed Laplacian can be defined \cite{cho} 
as $-\Delta = d\d + \d d$. 
In this form, the invariance under the isometries $U_h(sl(2))$ is obvious. 
To calculate it explicitely, we introduce 
\cite{cho_madore} derivations $e_a$ dual to the 1-forms 
$\theta^a$, defined by
$$
\begin{array}{ll}
e_1 x = y, &e_1 y = 0,\\[2pt]
e_2 x = 0, &e_2 y = -  y.
\end{array}
$$
In terms of them the Laplace operator $\Delta_h$ can be 
written  as
\be
- \Delta_h \phi = e_1^2 \phi + e_2^2 \phi + e_2 \phi, \label{Laplace}
\qquad \phi \in \c{A}_h.
\ee

First we recall the calculation of the propagator in the
commutative case. In the commutative limit $\Delta_h$ tends to the
ordinary Laplace operator on the Lobachevsky plane:
\be
\lim_{h\to 0} \Delta_h = \tilde \Delta = 
- \tilde y^2 (\partial_{\tilde x}^2 + \partial_{\tilde y}^2).
\ee
Here $(\tilde{x},\tilde{y})$ are the commutative limits
of the operators $(x,y)$. The spectrum of $\Delta_h$ in the commutative 
limit is given ~\cite{Ter85} by the eigenvalue equation
\be
\tilde\Delta \phi (\tilde{x}, \tilde{y}) = 
\lambda_{k,\kappa} \phi (\tilde{x}, \tilde{y}).                  \label{3.2}
\ee                                                 
By the separation of variables 
$\phi(\tilde{x}, \tilde{y}) = f(\tilde{x})g(\tilde{y})$ one finds the 
differential equations
\berr
&&\partial_{\tilde{x}}^2 f(\tilde{x}) = 
- k^2f(\tilde{x}),                                        \label{3.3}\\[4pt]
&&\tilde{y}^2 \partial_{\tilde{y}}^2 g(\tilde{y}) = 
(k^2 \tilde{y}^2 - \lambda_{k,\kappa}) g(\tilde{y}) \label{3.4}
\err
where $k \in \mf{R}$. We define 
$\kappa^2$ by 
$$
\lambda_{k,\kappa} = \kappa^2 + \frac 14. 
$$
The eigenvalues $\lambda_{k,\kappa}$ do not in fact depend on 
$k$ and are infinitely degenerate.  If we set then $z = ik\tilde{y}$ and
$g(\tilde{y}) = \sqrt{z}J(z)$, Equation~(\ref{3.4}) becomes the Bessel
equation
\be
J''(z) + \frac{1}{z}J'(z) + (1 + \frac{\kappa^2}{z^2})J(z) = 0.   \label{3.5}
\ee                                                      
A normalized set of eigenfunctions for the Laplace operator is given by
\be
\phi_{k,\kappa}(\tilde{x}, \tilde{y}) = 
e^{ik\tilde{x}} \pi^{-3/2}
\sqrt{\kappa \sinh \pi \kappa}  
\sqrt{\tilde{y}}K_{i\kappa}(|k|\tilde{y})                       \label{3.6}
\ee  
with $\kappa > 0$ and $k \neq 0$.  The case $\kappa < 0$ can be excluded
since
$$
K_{-\nu }(|k|\tilde{y}) = K_{\nu}(|k|\tilde{y}).
$$
The case $k = 0$ is also excluded since when $\tilde{y} \to 0$
\be
K_{i\kappa }(| k | \tilde{y}) \to \frac{1}{2}
\Gamma(i\kappa )\left(\frac{2}{|k| \tilde{y}}\right)^{i\kappa } + 
\frac{1}{2} \Gamma (-i\kappa )
\left(\frac{2}{|k| \tilde{y}}\right)^{-i\kappa}.              \label{expan}
\ee
If we set $\tilde x^i = (\tilde x, \tilde y)$ the completeness relation 
can be written as
\be
\delta^{(2)} (\tilde x^i - \tilde x^{i\prime}) = \int_{-\infty}^{+\infty} 
\int_0^\infty  \phi_{k,\kappa}(\tilde x, \tilde y)
\phi^*_{k,\kappa}(\tilde x^\prime, \tilde y^\prime) dk d\kappa
\ee
and the propagator of $(\Delta + \mu^2)$ is given by
\be
G (\tilde x^i, \tilde x^{i\prime}) = 
\int_{-\infty}^{+\infty} 
\int_0^\infty  \frac{\phi_{k,\kappa}(\tilde x, \tilde y)
\phi^*_{k,\kappa}(\tilde x^\prime, \tilde y^\prime)} 
{\kappa^2 + \frac 14 + \mu^2} dk d\kappa.  
                                                             \label{Green_cl}
\ee
Consider now the noncommutative case.  Notice that although
the classical Lobachevsky plane is invariant under the reflection 
$\tilde x \to - \tilde x$, this is no longer the case when $h \neq 0$.
By ordering again any monomial in $\A$ in the form $\phi(x|y)$,
one can formally separate the variables in the
eigenvalue problem as before \cite{cho} and the eigenvalue equation can be
decomposed into two differential equations. The equations for the factor
$f(x)$ are given by
\be
\begin{array}{l}
e_1^2 f(x) = - L_+^2y^2 f(x),                          \\[4pt]
e_1^2 f(x) = - L_-^2 f(x) y^2
\end{array}                                                 \label{3.7}
\ee
where $L_\pm \in \mf{R}$.  Since the commutation relations $[y, e_2]$ and
$[\tilde{y}, \tilde y \partial_{\tilde{y}}]$ are of the same form, the
differential equation for $g(y)$ has the same form as that of
(\ref{3.4}) even though the algebra has changed:
$$
(e_2^2 + e_2) g(y) = (L_\pm^2 y^2 - \lambda_{k,\kappa}) g(y).
$$
Consider the functions
$$
L_{\pm}(k) = \frac{e^{\pm 2hk} - 1}{2h}.
$$
For any
$k \in \mf{R}$ let $e^{ikx}$ be defined as a formal power series in the
element $x$. 
Then from the action of $e_1$ on $x$ it follows that
\be
e_1 e^{ikx} = i L_+(k) y e^{ikx} =  i L_-(k) e^{ikx} y.     \label{3.9}
\ee    
where we have used 
\be
e^{ikx} f(y) = f(e^{2hk}y) e^{ikx}.                           \label{com-rel}
\ee                  
The solution of Equation~(\ref{3.7}) is given therefore by
\be
f(x) = e^{ikx},   \qquad L_\pm = L_{\pm}(k).                 \label{3.10}
\ee                                                       
A family of formal solutions of the eigenvalue equation on the
quantum Lobachevsky plane which tend to normalized functions in the
commutative limit is given for $k \neq 0$, $\kappa >0$ by
\be
\phi_{k,\kappa}(x,y) = 
\pi^{-3/2} \sqrt{\kappa \sinh \pi \kappa}  
e^{ikx}\sqrt{y}K_{i\kappa}(|L_-(k)|y).                      \label{eigenfunct}
\ee  
Thus $L_-(k)$ plays the role of the linear momentum 
associated to $x$. Although $|k|$
remains invariant under the map $k \to -k$ this is not the case for
$|L_-(k)|$, a fact which is a manifestation of the breaking of parity by the
commutation relations. Moreover, the range of the momentum $L_-(k)$
in $x$ direction appears at this point
to be limited to the region $(-\frac 1{2h},\infty)$.
We will come back to this in a moment.

Define the 1--particle Hilbert space $\H^{(0)}$ to be 
generated by the (improper) basis 
$$
\phi_{k,\kappa}(x|y) =
\pi^{-3/2} \sqrt{\kappa \sinh \pi \kappa}  
e^{ikx}\sqrt{y}K_{i\kappa}(|L_-(k)|y),
$$ 
for $k \in (-\infty, \infty)$
and $\kappa >0$. The inner product on this space  
should be invariant under $U_h(sl(2))$, which means that the star--structure
(\ref{h_star}) of $U_h(sl(2))$ is induced by the adjoint
of operators on the Hilbert space, as classically.
As explained in the previous section, one such
inner product is given by 
\be
\<f(x|y), g(x|y) \>_{(0)} = 
     \int :f(x|y)^* g(x|y): d\mu,
\ee
where the latter is the classical integral after normal ordering,
i.e. $x$ should be commuted to the left of $y$ before taking the integral. 
It is clear that the Laplacian is (formally) a symmetric operator, since
$\d$ is the (formal) adjoint of $d$.

Now we can calculate the norm of the eigenstates as 
\berr
\<\phi_{k,\kappa}, \phi_{k',\kappa'}\>_{(0)} &=& 
      \int :\phi_{k,\kappa}(x,y)^* \phi_{k',\kappa'}(x,y): d\mu \nonumber\\
      &=& \pi^{-3} \int \sqrt{\kappa\sinh(\pi\kappa)}
                               \sqrt{\kappa'\sinh(\pi\kappa')} \nonumber \\
      &&  \quad  :\sqrt{y} K^*_{i\kappa}(|L_-(k)|y) e^{-i(k-k')x}
          \sqrt{y} K_{i\kappa'}(|L_-(k')|y):  d\mu \nonumber \\ 
      &=& \d(k-k') \pi^{-3/2} \int_0^{\infty} dy \frac 1{y} 
       \sqrt{\kappa\sinh(\pi\kappa)}\sqrt{\kappa'\sinh(\pi\kappa')}\nonumber\\
      &&  \quad   K^*_{i\kappa}(|L_-(k)|y) K_{i\kappa'}(|L_-(k')|y) \nonumber\\
      &=& \d(k-k') \d(\kappa-\kappa'),
\err
as classically. The Hilbert space $\H^{(0)}$ can now be defined 
as the closure of normalizable
wave--packets build from this ``basis'' of eigenfunctions, 
which obviously define an isometry with the usual, 
undeformed Hilbert space of square integrable functions.

Using 2 braided copies of $\A$ as in the previous section, 
the propagator can be written as
\berr
G(x^i,x^{i\prime}) &=& \int_{-\infty}^{+\infty} 
          \int_0^\infty \; \frac{\phi_{k,\kappa}(x,y) 
         \phi^*_{k,\kappa}(x^\prime,y^\prime)}{\lambda_{k,\kappa} + \mu^2}
          dk d\kappa                          \nonumber \\
    &=& \pi^{-3} \int_{-\infty}^{+\infty} 
   \int_0^\infty \;( \lambda_{k,\kappa}+ \mu^2)^{-1}\kappa\sinh(\pi\kappa) 
                                                       \nonumber\\  
    && \quad   e^{ikx} \sqrt{y}K_{i\kappa}(|L_-(k)|y)
       \sqrt{y'} K^*_{i\kappa}(|L_-(k)|y') e^{-ikx'} dk d\kappa \\
    &=& \pi^{-3} \int_{-\infty}^{+\infty} 
   \int_0^\infty \; (\lambda_{k,\kappa}+ \mu^2)^{-1}\kappa\sinh(\pi\kappa)
                                                       \nonumber\\   
    && \quad  e^{ikx} e^{-ikx'}\sqrt{y}K_{i\kappa}(|L_+(k)|y)
       \sqrt{y'} K^*_{i\kappa}(|L_+(k)|y') e^{2hk} dk d\kappa. \label{green_0}
\err
We have used here (\ref{com-rel}), the identity 
$|L_-(k)| e^{2hk} = |L_+(k)|$ 
and the fact that the commutation relations (\ref{braid_alg}) 
between $y$ and $x'$ are the same as those between $y'$ and $x'$.

As is shown in  appendix A, the commutation relations (\ref{braid_alg})
between $x$ and $x'$ imply the following identity:
\be
e^{ikx} e^{-ikx'} = e^{iL_+(k) \d x},                      \label{lemma}
\ee
where we recall $\d x = x-x'$.
Together with (\ref{com-rel}), it follows
\berr
G(x^i,x^{i\prime}) &=& \pi^{-3} \int_{-\infty}^{+\infty} 
   \int_0^\infty \; (\lambda_{k,\kappa}+\mu^2)^{-1}\kappa\sinh(\pi\kappa) 
                                                       \nonumber\\  
    && \quad  e^{iL_+(k)\d x} \sqrt{y} K_{i\kappa}(|L_+(k)|y)
       \sqrt{y'} K^*_{i\kappa}(|L_+(k)|y') e^{2hk} dk d\kappa . 
\err
Now recall that the subalgebra generated by $\obar{y}, \d x, \d y$
is abelian. Thus we can treat it as ordinary function algebra, and
change variables to $p=L_+(k), \quad dp = e^{2hk} dk$. We have then
\be
G(x^i,x^{i\prime}) = \pi^{-3} \int_{-\frac 1{2h}}^{+\infty} dp 
 \int_0^\infty d\kappa \;(\lambda_{p,\kappa}+\mu^2)^{-1}\kappa\sinh(\pi\kappa)
  e^{ip\d x} \sqrt{y} K_{i\kappa}(|p|y) \sqrt{y'} K^*_{i\kappa}(|p|y').
\ee
Recall that $\lambda_{p,\kappa}$ does not depend on $p$. 
The integrand is therefore exactly the same as classically (See 
(\ref{Green_cl}) and (\ref{3.6})); only the integration limit of $p$ has 
changed. 

This is actually a rather strange result. Since the Laplacian is invariant
under $U_h(sl(2))$, one should expect that the propagator be also 
invariant under this algebra (this is made more
explicit in appendix B), which implies that it is a function of the
invariant distance, which is the same as classically as we have seen. 
We just found that it is ``almost'', but not quite: 
if the integration limits were the same as classically, it would of course be
invariant, but the integration bound $-\frac 1{2h}$ for $p$  
spoils invariance. How is this possible?

The only explanation seems to be that the representation of 
$U_h(sl(2))$ on the Hilbert space generated by the 
eigenfunctions (\ref{eigenfunct}) is not a $\ast$--representation, 
i.e. the generators of $U_h(sl(2))$ are not represented as 
(anti)self--adjoint operators. 
A similar phenomenon is known to happen in the case of the $q$--deformed
quantum line \cite{schwenk}, where one has to consider reducible 
Hilbert space representations in order to 
obtain self--adjoint representations of the quantum algebra. 

In fact, we can find such a ``extended'' Hilbert space here as well. Let 
\berr
\phi_{k,\kappa}^{(n)}(x|y) &=& \pi^{-3/2} \sqrt{\kappa \sinh \pi \kappa}  
   e^{ikx}\sqrt{y}K_{i\kappa}(|L_-(k)|y)    \quad {\rm for}  \nn\\
 && k \in (-\infty, \infty)+\frac {i n\pi}{2h} \quad {\rm and} \quad \kappa >0,
   \label{basis_n}
\err
and let $\H^{(n)}$ be the closure of normalizable wavepackets 
built from (\ref{basis_n}),  with an inner product defined as
\be
\<f(x|y), g(x|y) \>_{(n)} :=\int :f(x,y)^* e^{n\pi x/h} g(x,y): d\mu.
\ee
It is easy to see that this inner product is 
invariant under the sub--Hopf algebra of $U_h(sl(2))$ generated by 
$\{G^2, J^3G,J^-G\}$ which commutes with 
$e^{n\pi x/h}$, and that the Laplacian is still symmetric since 
$d e^{n\pi x/h} =0$. 

As Hilbert space, $\H^{(n)}$ is of course equivalent to $\H^{(0)}$, but not 
as representation of $U_h(sl(2))$. For example, consider the above
``plane--wave'' states in $\H^{(1)}$: they are the eigenstates of 
the Laplacian with momentum $L_-(k)$ in $x$ direction 
in the interval $(-\infty, -\frac 1{2h})$, which were ``missing'' above.
We can calculate the inner product on $\H^{(1)}$:
\berr
\<\phi_{k,\kappa}^{(1)}, \phi_{k',\kappa'}^{(1)}\>_{(1)} &=& 
      \int :\phi_{k,\kappa}^{(1)}(x,y)^* e^{\frac{n\pi x}h} 
                            \phi_{k',\kappa'}^{(1)}(x,y): d\mu \nonumber\\
      &=&  \pi^{-3} \int \sqrt{\kappa\sinh(\pi\kappa)}
                               \sqrt{\kappa'\sinh(\pi\kappa')} \nonumber \\
      &&  \quad  :\sqrt{y} K^*_{i\kappa}(|L_-(k)|y) e^{-i(k-k')x}
          \sqrt{y} K_{i\kappa'}(|L_-(k')|y):  d\mu \nonumber \\ 
      &=& \d(k-k') \pi^{-3/2} \int_0^{\infty} dy \frac 1{y} 
       \sqrt{\kappa\sinh(\pi\kappa)}\sqrt{\kappa'\sinh(\pi\kappa')}\nonumber\\
      &&  \quad   K^*_{i\kappa}(|L_-(k)|y) K_{i\kappa'}(|L_-(k')|y) \nonumber\\
      &=& \d(k-k') \d(\kappa-\kappa').
\err 
One can now repeat the calculation (\ref{green_0}) for $\H^{(1)}$,
\berr
G^{(1)}(x^i,x^{i\prime}) &=& \int dk 
         \int_0^\infty d\kappa 
      \; \frac{\phi_{k,\kappa}^{(1)}(x,y) {\phi^{(1)}}^*_{k,\kappa}
        (x^\prime,y^\prime)  e^{\frac{n\pi x}h}}{\lambda_{k,\kappa}+\mu^2} 
                \nonumber \\
      &=& \pi^{-3} \int dk \int_0^\infty d\kappa  
             \; (\lambda_{k,\kappa}+\mu^2)^{-1}\kappa\sinh(\pi\kappa) 
                 \cdot \nonumber\\  
    && \quad \cdot  e^{iL_+(k)\d x} \sqrt{y} K_{i\kappa}(|L_+(k)|y)
       \sqrt{y'} K^*_{i\kappa}(|L_+(k)|y') e^{2hk} . 
\err
Changing again  variables to $p=L_+(k), \quad dp = e^{2hk} dk$, we obtain
\be
G^{(1)}(x^i,x^{i\prime}) = \pi^{-3} \int_{-\infty}^{-\frac 1{2h}} dp 
 \int_0^\infty d\kappa \;(\lambda_{p,\kappa}+\mu^2)^{-1}\kappa\sinh(\pi\kappa)
  e^{ip\d x} \sqrt{y} K_{i\kappa}(|p|y) \sqrt{y'} K^*_{i\kappa}(|p|y').
\ee
This is precisely the missing piece in order to obtain an invariant propagator.
We therefore define the ``extended'' Hilbert space 
to be the direct orthogonal sum
\be
\c{H}:=\c{H}^{(0)} \oplus  \c{H}^{(1)}.
\ee
Then on $\c{H}$, the propagator is invariant and exactly as classically,
\be
G^{\c{H}}(x^i,x^{i\prime}) = \pi^{-3} \int_{-\infty}^{\infty} dp 
  \int_0^\infty d\kappa \; (\lambda_{p,\kappa}+\mu^2)^{-1}\; 
  \kappa\sinh(\pi\kappa)  
  e^{ip\d x} \sqrt{y} K_{i\kappa}(|p|y) \sqrt{y'} K^*_{i\kappa}(|p|y').
\ee

\sect{Discussion}

It was found in \cite{madore_etc} that the propagator on the $h$--deformed 
Lobachevsky plane is finite if one uses the usual, ``unbraided'' tensor
product; similarly for  a noncommutative flat plane. However, this tensor
product ``breaks'' the invariance under the quantum group. 

In this paper, we have first seen that the $h$--deformation is not
a trivial deformation, even though it is just a twist of the 
undeformed case; it turns out that the structure of the 
Hilbert space is modified.
If this is done properly and the covariant, 
braided tensor product is used,
then the propagator turns out to be the same as classically; in
particular  it is divergent. 
This is certainly disappointing, since the main reason for considering the
deformed manifolds is to ``smear'' the points, thereby regularizing the 
UV divergences. 
It is not entirely clear how to understand this result.
It could be that the algebra is simply not noncommutative enough. 
Another interpretation 
might be that the identification of the ``distances'' $\d x, \d y$ in the
braided tensor product is not satisfactory, in particular since they do  not
commute with the ``average'' values $\obar{x}, \obar{y}$. See also the related
discussion in \cite{madore_etc}. 
In other words, the nontrivial (braided) commutation relations between
different copies of the quantum space which are required by the
covariance under a quantum group imply some kind of interaction. 
It seems that the physical meaning of 
the braided tensor product is not completely understood, and deserves
further investigation. 

\sect{Acknowledgements}
The authors would like to thank Prof. P. Kulish for interesting
discussions. This work was
partially supported by the DAAD under the PROCOPE grant number
PKZ~9822848.

\sect{Appendix A}
We prove that $e^{ikx} e^{-ikx'} = e^{iL_+(k)\d x}$, where 
$\[x,x'\]  = 2ih (x-x')$. 

From $[x,\d x] = -2ih \d x$ it follows 
$[x, ... [x,x'] ... ] = -(-2ih)^n \d x$ for $n$ commutators, and thus
\be
e^{ikx} x' e^{-ikx} = x'- (e^{2hk}-1) \d x.          \label{l_1}
\ee
Let $f_k(x,x') = e^{ikx} e^{-ikx'}$. Using (\ref{l_1}), we find
\berr
\frac d{dk} f_k(x,x') &=& ix e^{ikx} e^{-ikx'} - e^{ikx} ix' e^{-ikx'} \nn\\
           &=& i(x-x'+ (e^{2hk}-1)\d x)f_k(x,x') \nn\\
           &=& e^{2hk} i\d x f_k(x,x').           \label{l_2}
\err
Consider $\frac{dp(k)}{dk} = e^{2hk}$, with the solution
$p(k) = \frac 1{2h} (e^{2hk}-1)$. Then (\ref{l_2}) is
equivalent to 
$$
\frac{df_p(x,x')}{dp} = i\d x f_p(x,x'),
$$
with the solution $f_p(x,x') = e^{ip \d x}$. Therefore
\be
e^{ ikx } e^{-ikx'} = e^{\frac i{2h} (e^{2hk}-1)\d x} = 
                  e^{iL_+(k) \d x}.
\ee

\sect{Appendix B}
We explain without mathematical rigour why the propagator should be
invariant under $U_h(sl(2))$ if the latter is implemented via 
a $\ast$--representation. We assume the spectral decomposition 
$\one = \int_n \phi_n \tens \phi_n^*$ where $\Delta \phi_n  =\l_n \phi_n$.

For $u \in U_h(sl(2))$, one has 
\berr
u\cdot G(z,z') &=& u\cdot \(\int_n \l_n^{-1} \phi_n(z) 
                                           \phi_n^*(z')\) \nonumber\\
  &=& \int_n \l_n^{-1} (u_{(1)}\cdot \phi_n(z))(u_{(2)}\cdot \phi_n^*(z'))
\err
where $z$ stands for $(x,y)$.
Let 
$$
u\cdot \phi_n = \int_l \phi_l \pi_{ln}(u).
$$
Since $\pi$ is a $\ast$--representation one has 
$\pi_{ln}(u^*) = \pi_{nl}^*(u)$. We claim that this implies that 
\be
u\cdot \phi_k^* = \int_l \pi_{kl}(Su) \phi_l^*
\ee
where $Su$ is the antipode of $u$. Indeed using $U_h(sl(2)) \smash \A$, 
one has
\berr
u\cdot \phi_k^* &=& u_{(1)} \phi_k^* Su_{(2)} = \(S^{-1}(u_{(2)}^*) 
         \phi_k u_{(1)}^*\)^* = \((S^{-1}u^*)\cdot \phi_k\)^* \nonumber\\
   &=&  \(\int_l \phi_l \pi_{lk}(S^{-1}u^*)\)^* = 
              \int_l \phi_l^* \pi_{lk}^*(S^{-1}u^*) 
         = \int_l \pi_{kl}(Su) \phi_l^*
\err
where the $\ast$--representation property was used in the last line, 
as well as standard identities for $\ast$--Hopf algebras. 
Now we can conclude that
\berr
u\cdot G(z,z') &=&  \int_{n,l,k} \l_n^{-1} \phi_l(z) \pi_{ln}(u_{(1)}) 
                    \pi_{nk}(S u_{(2)}) \phi_k^*(z') \nonumber\\
  &=& \epsilon(u)  \int_{n,l,k}  \phi_l(z) \d_{lk} \phi_k^*(z') = 
      \epsilon(u) G(z,z').
\err

\end{document}